\documentclass[hidelinks,onefignum,onetabnum]{siamart220329}

\usepackage{bbm}
\usepackage{subcaption}
\usepackage{amssymb}
\usepackage{listings}

\newcommand{\e}{\mathrm{e}}

\newcommand{\R}{\mathbb{R}}
\newcommand{\C}{\mathbb{C}}
\newcommand{\N}{\mathbb{N}}
\newcommand{\bn}{\mathbf{n}}
\newcommand{\bx}{\mathbf{x}}

\newcommand{\btau}{\boldsymbol{\tau}}

\newcommand{\by}{\mathbf{y}}
\newcommand{\br}{\boldsymbol{r}}

\newcommand{\de}{\,\mathrm{d}}

\newcommand{\tvarphi}{\widetilde \varphi}
\newcommand{\real}{\operatorname{Re}}
\newcommand{\imag}{\operatorname{Im}}

\newsiamthm{remark}{Remark}

\graphicspath{{./figures/}}

\title{A complex-scaled boundary integral equation for time-harmonic water waves}

\author{\hspace{1cm}Anne-Sophie Bonnet-Ben Dhia\thanks{POEMS, CNRS, Inria, ENSTA Paris, Institut Polytechnique de Paris,  France.}
\and Luiz M. Faria\footnotemark[1]\and\newline
 Carlos P\'erez-Arancibia\thanks{Department of Applied Mathematics, University
of Twente, Enschede, The Netherlands}}

\date{\today}

\begin{document}

\maketitle

\begin{abstract}  
This paper presents a novel boundary integral equation (BIE) formulation for the two-dimensional time-harmonic water-waves problem. It utilizes a complex-scaled Laplace's free-space Green's function, resulting in a BIE posed on the infinite boundaries of the domain. The perfectly matched layer (PML) coordinate stretching that is used to render propagating waves exponentially decaying, allows for the effective truncation and discretization of the BIE unbounded domain. We show through a variety of numerical examples that, despite the logarithmic growth of the complex-scaled Laplace's free-space Green's function, the truncation errors are exponentially small with respect to the truncation length. Our formulation uses only simple function evaluations (e.g. complex logarithms and square roots), hence avoiding the need to compute the involved water-wave Green's function. Finally, we show that the proposed approach can also be used to find complex resonances through a \emph{linear} eigenvalue problem since the Green's function is frequency-independent.
  
\end{abstract}


\section{Introduction}

Wave propagation and scattering problems in physics and engineering often involve unbounded domains that include infinite boundaries. These types of models are particularly important in various scenarios, such as seismic waves interacting with the free surface, outdoor sound propagation, guided (acoustic, electromagnetic, elastic) wave phenomena, and water waves propagating at the ocean's free surface, to mention a few. This paper specifically focuses on the linear (finite-depth) time-harmonic water wave problem, which features two infinite boundaries modeling the bottom topography and the free surface. 
This classical problem can be traced back to the work of Cauchy~\cite{cauchy1816theory} and originated from the study of the hydrodynamics of floating and submerged bodies. Extending beyond water waves, the simple model considered in this paper could serve as a prototypical example for investigating a wide range of simplified wave problems where the oscillations of volumetric waves can be neglected at the scale of surface waves. This is typically the case in plasmonics, where the relevant physical phenomenon is the presence of (subwavelength) surface plasmon polaritons~\cite{maier2007plasmonics}.

While boundary integral equation (BIE) methods are well known for their effectiveness in handling unbounded domains, they encounter notable challenges when applied to problems involving infinite boundaries. BIEs naively formulated using free-space Green's functions over the entire domain inherently result in equations defined on the domain's infinite boundaries. Therefore, suitable truncation techniques become essential to enable their numerical solution, provided, of course,  sufficient decay of the integrands exists. The simple approach of abruptly truncating the computational domain at some distance far from the wave sources, offers limited success in the case of the Helmholtz equation due to the slow decay of the free-space Green's function, and becomes completely ineffective when applied to two-dimensional water wave problems due to the logarithmic growth of Laplace's free-space Green's function. (We mention in passing, that  a domain decomposition approach that uses the free-space Green's function together with a BIE formulation over a bounded domain, has been pursued in \cite{mei1978numerical,yeung1982numerical}. In those works, a Dirichlet-to-Neumann map based on a modal series expansion of the scattered velocity potential is used to effectively reduce the problem to a bounded domain, an idea that later seemed to have influenced the development of finite element methods for water wave problems~\cite{lenoir1988localized}.)

Leveraging the geometric simplicity of the infinite boundaries involved, many BIE formulations for water wave problems rely instead on the use of problem-specific Green's functions that exactly satisfy the relevant boundary conditions over infinite flat portions of the domain's boundaries~\cite{angell1991recent,angell1986integral,black1975wave,chakrabarti2001application,fenton1978wave,mackay2021green}. While these formulations yield BIEs posed on a bounded portion of the domain's boundary, the numerical evaluation of the associated Green's function often becomes the major bottleneck in the solution process since simple expressions for the free surface Green's function are in general not available. Intricate numerical procedures have been developed in order to compute the free surface Green's function (see e.g.~\cite{xie2018comparison} and~\cite{mackay2021green} for a comparison of  different techniques and recent developments on this subject). This approach typically entails Fourier-like integrals, non-standard special functions, and/or power series expansions, which in view of the large number of Green's function evaluations involved in the BIE solution process, can be computationally demanding~\cite{newman1985algorithms}. Furthermore, off-the-shelf acceleration techniques for the Laplace equation, such as the fast multipole method (FMM)~\cite{greengard1987fast}, cannot be directly applied to the free surface Green's function. Interestingly, an explicit closed-form expression (in terms of exponential integrals) for the infinite depth problem in two dimensions, is presented in~\cite{hein2010explicit}. Even in this relatively simple case, developing a multipole expansion is non-trivial~\cite{perez2012fast}. Unfortunately, as far as the authors are aware, there are no closed-form expressions (nor multipole expansions) available for the two-dimensional finite depth problem as well as for the three-dimensional finite and infinite depth problems. We refer the reader to~\cite{martin2006multiple,hoernig2010green,linton2001mathematical} where these topics are discussed in great detail. In addition to the aforementioned efficiency issues, it is important to mention that problem-specific Green's functions always incorporate image sources which are essential for enforcing the free surface and Neumann boundary conditions at flat boundaries. Their presence, however, imposes significant constraints on the types of domains to which the associated BIE formulation can be applied. In fact, domains suitable for the application of such BIE formulations must be strictly confined within the strip defined by the flat free surface and the flat bottom topography.

Over the last decade, there have been notable advancements in the development of novel BIE formulations for Helmholtz wave problems involving unbounded interfaces~\cite{bruno2016windowed,lu2018perfectly}. These formulations revisit the idea of utilizing the simple efficient-to-evaluate free-space Green's function and truncating the resulting BIE posed on infinite boundaries, but they incorporate mechanisms to effectively account for the correct ``outgoing" solution behavior outside the bounded region of interest. In detail,~\cite{bruno2016windowed} put forth the so-called Windowed Green's function (WGF) method, which features a real-valued suitably-scaled smooth window function that multiplies the free-space Green's function kernel to effectively reduce the BIE to a bounded domain. In essence, this method exploits the fact that windowed integrands exhibiting decaying and oscillatory behavior can be neglected from the tails of the integrals over the infinite boundaries, introducing errors that decay super-algebraically fast as the size of the truncated region increases. Around the same time, a technique was introduced by Lu et al. in~\cite{lu2018perfectly} that utilizes a BIE formulation of the PDE boundary value problem transformed by a perfectly matched layer (PML) complex stretching and subsequently truncated to a bounded domain. As the WGF method, the resulting PML-BIE formulation is posed on a bounded portion of the infinite boundary and features a modified free-space Green's function that encodes the outgoing oscillatory character of the PDE solution. 

Although the WGF proves effective at solving Helmholtz problems in the context of layered and periodic media as well as waveguides~\cite{bruno2016windowed,bruno2017windowed2,bruno2017windowed,strauszer2023windowed}, its direct application to the water wave problem under consideration is not feasible. While it may be tempting to attribute the failure of the WGF method to the logarithmic growth of the free space Green's function, an argument can be made that by combining an easy-to-evaluate harmonic function with the free-space Green's function, a decaying Green's function can be formed (in a manner akin to the construction used in~\cite{chandler1998uniqueness}). However, even in that scenario, it is easy to see that the WGF method would fail to approximate the correct radiative water wave problem solution.  Indeed, when considering real boundary data, the WGF method produces real-valued (stationary-wave) solutions that fail to approximate the correct radiative water-wave problem solutions. We suspect that the failure of the WGF method in the context of the water waves problem is due to the fact that, unlike in the Helmholtz case, the Laplace free space Green's function contains no information about the ``correct`` radiation condition of the problem.

Building upon the PML-BIE method for the Helmholtz equation~\cite{lu2018perfectly}, this paper proposes a complex-scaled BIE method for the time-harmonic water wave problem in two dimensions. Our BIE formulation is derived from a direct integral representation of the analytically extended scattered velocity potential, obtained by suitably stretching the unbounded domain in the complex plane. This representation involves integrals over the infinite boundaries corresponding to the free surface and the bottom, and is obtained by applying  Green's third identity to the velocity potential and the logarithmic free-space Green's function over the complex stretched domain. Interestingly, despite the unboundedness of Laplace's Green's function along the complex stretched infinite boundaries, the resulting integral representation formula remains meaningful due to the provable exponential decay of the complexified velocity potential. Indeed, this decay is achieved by devising complex stretching coordinates on which the modal decomposition of the far-field scattered velocity potential, which encompasses a (constant-amplitude) propagating surface-wave mode and infinitely many evanescent modes, decays exponentially fast. Upon enforcing the boundary conditions and truncating the infinite boundaries, we then obtain a BIE that we can solve numerically via a high-order Nystr\"om method which suitably handles the complexified kernels. Through our numerical experiments, we demonstrate that the error in the computed solutions using this methodology decays exponentially fast as the length of the PML region increases. 

A noteworthy difference with the PML-BIE technique applied to the Helmholtz equation~\cite{lu2018perfectly} stems from the fact that, due to the lack of wave-like behavior in the fundamental solution of Laplace's equation, its complex extension still exhibits a logarithmic growth at infinity (a more detailed discussion is given in~\Cref{sec:complex-scaled-formulation}). This makes it not immediately obvious that truncation errors stemming from the method will converge to zero as the length of the PML region is increased. While for the exact solution it is straightforward to show that the truncation of the involved integral operators to a bounded domain yields exponentially small errors (for the analytic extension of the exact solution decays exponentially at infinity), proving a stability result on the proposed truncation technique remains an open question. The main difficulty lies in the fact that the kernels in the proposed integral equation do not decay exponentially at infinity (in fact they may even grow logarithmically in some cases), and thus the integral operators over the unbounded curves are only well-defined for densities decaying sufficiently fast. The numerical results of~\Cref{sec:numerical-results} suggest, however, that the lack of decay in the Green's function is not a fundamental problem in the method, and much like in the Helmholtz case we observe truncation errors that decay exponentially fast as the length of the PML layer is increased.

To further showcase the capabilities and possible advantages of the proposed methodology, we numerically investigate the challenging problem of the scattering of a surface wave by an infinite step bottom topography, where a direct application of a BIE formulation that makes use of the constant depth water-wave Green's function (satisfying free surface and Neumann boundary conditions at the flat infinite boundaries), is not feasible~\cite{athanassoulis1999consistent,belibassakis2004three,porter2000water}. Finally, we apply our methodology to approximate the resonant frequencies that arise in floating body problems~\cite{Haz-Len-1993,mciver1996example}. In contrast to alternative approaches for computing complex resonances based on the use of the water-waves Green's function, which yields a (non-linear) eigenvalue problem where the frequency appears in the kernels of the integral operators~\cite{Haz-Len-1993}, our PML-transformed formulation gives rise to a linear (generalized) eigenvalue problem that can be efficiently solved using standard matrix eigenvalue solvers.

The paper is organized as follows. In~\Cref{sec:prob_setup} we provide the setup of the problem under consideration, recalling some standard results of water-waves theory. \Cref{sec:complex-scaled-formulation} focuses on the construction of the complex stretching based on the modal decomposition of the scattered velocity potential, and the associated complex-scaled BIE formulation. In~\Cref{sec:Nyst_meth}, we introduce the Nystr\"om method employed in our approach. Next, \Cref{sec:numerical-results} presents a variety of numerical experiments designed to validate the proposed methodology and also demonstrate its practical applicability (the method's implementation, as well as the scripts used to generate the figures, are publicly available in an accompanying GitHub repository\footnote{\href{https://github.com/maltezfaria/HarmonicWaterWaves/tree/v0.1}{https://github.com/maltezfaria/HarmonicWaterWaves/tree/v0.1}}; all of the numerical results of~\Cref{sec:numerical-results} can be reproduced by following the instructions there). In~\Cref{sec:open-problems}, we briefly discuss some open problems and theoretical difficulties. Finally, we present in~\Cref{sec:conclusions} some concluding remarks and possible future directions.

\section{Problem setup}\label{sec:prob_setup}

This work focuses on linear and time-harmonic water waves in a bath of finite depth $d$. For the sake of presentation simplicity, in what follows we derive dimensionless equations, using the depth $d$ as a characteristic length scale (this means, in particular, that the dimensionless depth is one). 

\begin{figure}
\includegraphics[scale=1]{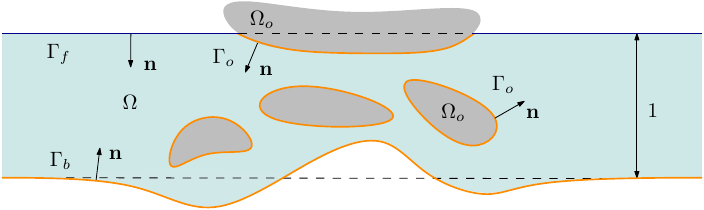}    
  \centering	  
  \caption{Schematics of the domain and main variables involved in the setup of water wave problem~\eqref{eq:water-waves-system} under consideration. \label{fig:schematic-geometry}}
\end{figure}

Throughout this paper, we let $\Omega_f \subset \mathbb{R}^2$ denote the infinite fluid domain bounded above by the free surface $\Gamma_s$ and below by a bottom topography $\Gamma_b$. We also let $\Omega_o \subset \mathbb{R}^2$
denote the domain occupied by bounded solid obstacles (either fully or partially) immersed in the fluid, with the immersed part of their boundary denoted by $\Gamma_o:=\partial \Omega_o\cap\overline\Omega_f$. More precisely, the unperturbed free surface is assumed to be placed at $\R\times\{0\}$ so that $\Gamma_f := \R \times \{0\} \setminus
\overline\Omega_o$.  The bottom topography, on the other hand, is assumed to be given by a constant depth $1$ outside of a bounded domain, i.e., there exists an $L_0>0$ such that $\{|x_1|>L_0\}\times \{-1\}\subset\Gamma_b$. In the derivations that follow we use $L_1>\max\{L_0,\operatorname{diam}\Omega_o\}$. %
A schematic of the geometry is shown in~\Cref{fig:schematic-geometry}.

Here we seek to solve the following boundary value problem:
\begin{subequations}
  \label{eq:water-waves-system}
  \begin{align}  
    \Delta \varphi(\bx) &=0,\phantom{(\bx)} \quad\ \bx:=(x_1,x_2) \in \Omega:=\Omega_f\setminus\overline\Omega_o,\\
    \nabla \varphi(\bx) \cdot \bn(\bx) + \nu \varphi(\bx)&=f_1(\bx), \quad \bx \in \Gamma_f:=\R \times \{0\} \setminus
    \overline\Omega_o,\label{eq:impedance_bc}\\
    \nabla \varphi(\bx) \cdot \bn(\bx) &=f_2(\bx), \quad \bx \in \Gamma_b \cup \Gamma_o,\\
    \label{eq:radiation-condition}  
    \lim_{x_1 \to \pm \infty} \partial_{x_1} \varphi(\bx) \mp i k \varphi(\bx)  &= 0,
  \end{align}
\end{subequations}
where letting $\omega>0$ denote the angular frequency and $g$ the
gravitational acceleration, the term $\nu = (\omega^2d)/g$ in the Robin boundary condition~\eqref{eq:impedance_bc} is a dimensionless parameter, while the constant $k > 0$ in the radiation condition~\cref{eq:radiation-condition} is uniquely determined by the dispersion relation
\begin{align}
\label{eq:dispersion-relation}
k\tanh(k) = \nu.
\end{align}
As usual, $\bn\in\R^2$ denotes the unit normal which points toward the interior of the
fluid domain $\Omega$ (see~\Cref{fig:schematic-geometry}).  It is worth noting that the unknown velocity potential $\varphi$ is  related to
the fluid velocity $\boldsymbol{u}$ through $\boldsymbol{u}(\bx,t) = \operatorname{Re}
\left(\e^{-i\omega t} \nabla \varphi \right)$, where~$t$ is the time variable.

We here assume that the boundary data $f_1$ and $f_2$ are (compactly) supported in $\Gamma_f$ and $\Gamma_b$ vanishing over $\{|x_1|>L_2\}\times\{0\}$ and $\{|x_1|>L_2\}\times\{-1\}$, respectively. Note that the step-like bottom topography problem considered in~\Cref{sec:numerical-results} is not covered by the previous description; nevertheless, a slight modification of the assumptions made here allows us to extend the proposed PML-BIE approach to this challenging problem.

In order to construct a suitable complex scaling for problem~\eqref{eq:water-waves-system}, we rely on the analytic extension of the PDE solution~$\varphi$. To justify the extension and gain insight into the construction of the complex stretching,  in what follows in this section we recall some classical results about the modal decomposition of $\varphi$ within semi-infinite strips of constant unit depth~\cite{kuznetsov2002linear,linton2001handbook}. Indeed, letting $L = \max\{L_1,L_2\}$, the exact solution $\varphi$ of~\eqref{eq:water-waves-system} admits a modal decomposition within the domain $\{|x_1|>L\}\times(-1,0)$. To define the modal decomposition we first introduce the sequence $\left\{\gamma_n\right\}_{n\geq 1}$, such that
\begin{equation}
  \label{eq:disp-evanescent}
  \gamma_n\tan(\gamma_n)=-\nu,\quad n\pi-\frac{\pi}{2}<\gamma_n <n\pi.
\end{equation}
Then one can prove that, with a suitable choice of the normalization constants
$(a_n)$, the functions $\{ u_n \}_{n\geq 0}$ defined by 
\begin{align}
\label{eq:transversemodes}
u_0(x_2)&=a_0\cosh\left(k(x_2+1) \right),&\\
u_n(x_2)&=a_n\cos(\gamma_n (x_2+1)), \quad n\geq 1,
\end{align}
with $k>0$ being the solution of~\cref{eq:dispersion-relation} and $\gamma_n$ the solution of~\cref{eq:disp-evanescent}, form an orthonormal basis of $L^2([-1,0])$. Invoking the radiation condition~\cref{eq:radiation-condition} to select the correct behavior at $\pm \infty$, the following modal expansion for the solution $\varphi$ of~\eqref{eq:water-waves-system} can be readily obtained:
\begin{align}
  \label{eq:modaldecomp}
  \varphi(\bx) = \begin{cases}
  \displaystyle A_0^+u_0(x_2)\e^{ik(x_1-L)} + \sum_{n\geq 1}A_n^+u_n(x_2)\e^{-\gamma_n(x_1-L)}, &\bx\in (L,\infty)\times [-1,0],\medskip\\
  \displaystyle A_0^-u_0(x_2)\e^{-ik(x_1+L)} +\sum_{n\geq 1}A_n^-u_n(x_2)\e^{\gamma_n(x_1+L)}, &\bx\in (-\infty,-L)\times [-1,0],
  \end{cases}
\end{align}
where 
\begin{align}\label{eq:modal-decomp-amplitude}
  A_n^\pm=\int_{-1}^0\varphi(\pm L,x_2)\overline{u_n(x_2)}\de x_2,\quad n\geq 0.  
\end{align}

The well-posedness of~\eqref{eq:water-waves-system} can be studied by reducing the problem to a bounded domain via a DtN map approach derived from the modal expansion~\eqref{eq:modaldecomp}. Indeed, it can be established that the problem~\eqref{eq:water-waves-system} is well posed for all frequencies except for a countable set (see e.g.~\cite{thomas1977eigenfunction}).

Clearly, on each connected component of the exterior domain $\{|x_1|>L\}\times[-1,0]$ the velocity potential~$\varphi$ is a superposition of the outgoing surface wave and of infinitely many evanescent modes. 

\section{Complex-scaled boundary integral equation}\label{sec:complex-scaled-formulation}

In this section, we derive a complex scaled integral equation for solving~\cref{eq:water-waves-system} with the desirable 
property that its solutions are exponentially decreasing as $|x_1| \to \infty$.
Using the modal decomposition~\cref{eq:modaldecomp}, we show
in Section~\ref{sec:complex-stretching} that the velocity potential $\varphi$ admits an analytic extension in $x_1$ into the complex plane for $|x_1| > L$ (i.e. outside the support of the perturbations). Using the analytic extension of the modal decomposition, we then derive in Section~\ref{sec:complex-scaled-bie} a modified water waves problem, as well as its boundary integral representation.

\subsection{Complex stretching}\label{sec:complex-stretching}

The modal decomposition~\cref{eq:modaldecomp} can be used to study analytic
properties of $\varphi$ with respect to the coordinate $x_1$. More precisely,
for $x_1>L$ (resp. $x_1<-L$) the solution $\varphi$
of~\cref{eq:water-waves-system} has an analytic extension to complex values of its first variable
with a real part larger than $L$ (resp. smaller than $-L$). To see why this is so, let us consider for instance the side $x_1>L$. It is
clear that for any $N\in\N$, the finite sum
\begin{align}
  A_0^+u_0(x_2)\e^{ik(z_1-L)}+\sum_{n= 1}^NA_n^+u_n(x_2)\e^{-\gamma_n(z_1-L)}
\end{align}
is well-defined for complex values of $z_1$, and is a holomorphic function of
$z_1$. Moreover, it converges when $N\to +\infty$,
uniformly in all compact subsets of $\real(z_1)>L$, to a function that we denote (with a slight abuse of notation) by $\varphi(z_1,x_2)$.  Finally, thanks to Morera's theorem, a uniform limit of a sequence of holomorphic functions is itself holomorphic, so we conclude that  $z_1\mapsto\varphi(z_1,x_2)$ is analytic within $\real(z_1)>L$. The exact same argument applies to the analytic extension $z_1\mapsto\varphi(z_1,x_2)$ for $\real(z_1)<-L$.

It follows then  by the analyticity of $\varphi$ that
\begin{align}  
  \label{eq:pml-laplace-equation}
  \left( \partial_{z_1}^2 + \partial_{x_2}^2 \right) \varphi(z_1,x_2) &=0, \quad (z_1,x_2) \in \{|\real z_1 | > L \} \times [-1,0].
\end{align}
We now introduce a (continuous and piece-wise differentiable) stretching path $\tau : \R \to \C$ along which the transformed function 
\begin{equation}\label{eq:transf_eq}
\tvarphi(x_1,x_2) := \varphi(\tau(x_1),x_2)    
\end{equation} decays exponentially fast as $|x_1| \to \infty$. In view of the modal decomposition~\cref{eq:modaldecomp}, the following conditions guarantee the exponential decay of $\tvarphi(\cdot,x_2)$:
\begin{subequations}
\label{eq:tau-conditions}
\begin{align}
    \imag\left(\tau(x_1)\right) &> 0 \quad \mbox{for} \quad x_1 > L\\
    \imag\left(\tau(x_1)\right) &< 0 \quad \mbox{for} \quad x_1 < -L
\end{align}        
\end{subequations}
Furthermore, in order to recover the original solution to~\cref{eq:water-waves-system} for $|x_1| < L$, we require
\begin{align}
    \tau(x_1) = x_1 \quad \mbox{for} \quad |x_1|<L. \tag{\ref{eq:tau-conditions}c} \nonumber
    \label{eq:tau-condition-identity}
\end{align}
Finally, we suppose that
\begin{align}
    0 < c < \real(\tau'(x_1)) < C < \infty \quad \mbox{for} \quad x_1 \in \R, \tag{\ref{eq:tau-conditions}d}
    \label{eq:tau-condition-real}
\end{align}
which implies that $\tau$ is injective. 

A simple application of the chain rule then yields the identity $\partial_{z_1} \varphi(\tau(x_1),x_2) = \partial_{x_1} \tvarphi(x_1,x_2)/\tau'(x_1)$, therefore, in view of~\eqref{eq:pml-laplace-equation},  $\tvarphi$ defined in~\eqref{eq:transf_eq} satisfies
\begin{subequations}
\begin{align}
  \label{eq:pml-laplace-real}
  \nabla \cdot (A(\bx) \nabla \tvarphi) &= 0  \quad \mbox{with} \quad A(\bx) :=   
  \begin{bmatrix}
    \tau'(x_1)^{-1} & 0 \\
    0 & \tau'(x_1)
  \end{bmatrix}, \quad \bx \in \Omega.
\end{align}
\label{eq:water-waves-system-pml}
\end{subequations}
(We have chosen a matrix formalism to simplify the discussion in the subsequent section.) Note that by~\cref{eq:tau-condition-identity}, $\tau' = 1$ for $|x_1| < L$,  so~\cref{eq:pml-laplace-real} becomes the Laplace equation in that case. Finally, we note that because $f_1$ and $f_2$, as well as the perturbations of the geometry vanish for $|x_1|>L$, it can be easily shown that $\tvarphi$ satisfies the same boundary conditions as $\varphi$, i.e.,
\begin{align}
    \label{eq:pml-free-surface-bc}
    \nabla \tvarphi(\bx) \cdot \bn(\bx) + \nu\tvarphi(\bx)&=f_1(\bx), \quad \bx \in \Gamma_f, \tag{\ref{eq:water-waves-system-pml}b}\\
    \label{eq:pml-solid-bc}  
    \nabla \tvarphi(\bx) \cdot \bn(\bx) &=f_2(\bx), \quad \bx \in \Gamma_b \cup \Gamma_o. \tag{\ref{eq:water-waves-system-pml}c}
\end{align}

It is worth mentioning at this point that the radiation condition~\cref{eq:radiation-condition} is simply replaced by the exponential decay of $\tvarphi$. Equations \cref{eq:water-waves-system-pml} comprise the PML-transformed water waves system that we seek to solve in lieu of the original system~\cref{eq:water-waves-system}.

\subsection{Boundary integral equation}\label{sec:complex-scaled-bie}
In order to recast the problem as a BIE by leveraging fundamental potential theory results, we first show that our PDE~\cref{eq:pml-laplace-real} is strongly elliptic in the sense of \cite{mclean2000strongly}. Indeed, given the form of the matrix function $A$ introduced in~\cref{eq:pml-laplace-real} and the stretching path property~\cref{eq:tau-condition-real}, it follows that
\begin{align}
    \real(\xi^* A \xi) = \real\left( \frac{1}{\tau'} |\xi_1|^2 + \tau' |\xi_2|^2 \right) > \min\left(c,\frac{1}{C}\right) |\xi|^2,\quad \xi \in \C^2,
\end{align}
which shows the desired property.

A crucial aspect of the subsequent BIE formulation is the  explicit relationship 
\begin{align}
  \label{eq:complex-Green-function}
\widetilde{G}(\bx,\by) := G(\btau(\bx),\btau(\by)) = -\frac{1}{4\pi}\log \left( (\tau(x_1) - \tau(y_1))^2 + (x_2 - y_2)^2 \right),
\end{align}
between the Green's function $\widetilde{G}$ associated with~\cref{eq:pml-laplace-real}, and the Laplace's fundamental solution $G$, along with the transformation $\tau$ (see~\Cref{sec:fundamental-solution}). Here, the principal logarithm is used with a branch cut on the negative real line and, for the sake of compactness of notation, we have introduced the vector-valued transformation 
$$\btau(\bx) := (\tau(x_1),x_2) : \R^2 \to \C \times \R.$$ %

With these definitions at hand, we are in position to derive the boundary integral representation for $\tvarphi$, defined in~\eqref{eq:transf_eq}, over the (unbounded) domain $\Omega$. We then begin by considering the truncated domain $\Omega_M := \left\{\bx \in \Omega : |x_1| <  M \right\}$, $M>0$. On $\Omega_M$, the following integral representation formula holds~\cite{mclean2000strongly}:
\begin{align}
  \label{eq:greens-representation-bounded}
  \tvarphi(\bx) = \mathcal{D}_{\partial \Omega_M}[\tvarphi](\bx) - \mathcal{S}_{\partial \Omega_M}[\nabla \tvarphi \cdot A \bn](\bx), \quad \bx \in \Omega_M,
\end{align}
where $\mathcal{S}_\Sigma$ and $\mathcal{D}_\Sigma$ denote respectively the single- and double-layer potentials, defined as
\begin{subequations}
\begin{align}
    \label{eq:SL-pot}
    \mathcal{S}_{\Sigma}[\sigma](\bx) &:= \int_\Sigma \widetilde{G}(\bx,\by) \sigma(\by) \de s(\by)\quad\text{and}\\
    \label{eq:DL-pot}
    \mathcal{D}_{\Sigma}[\sigma](\bx) &:= \int_\Sigma \left( \nabla_{\by} \widetilde{G}(\bx,\by) \cdot A(\by)\bn(\by) \right) \sigma(\by) \de s(\by),\quad\bx\in\R^2\setminus\Sigma.
\end{align}
\label{eq:potentials}
\end{subequations}
In~\cref{eq:potentials}, $\Sigma$ stands for a generic sufficiently regular (Lipschitz) curve. 

Using the fact that $\tvarphi$ decays exponentially fast as $|x_1| \to \infty$, while the integral kernels in~\cref{eq:greens-representation} grow at most logarithmically, the contribution from the lateral curves $\{ \pm M\} \times
[-1,0]\subset\partial\Omega_M$ in~\eqref{eq:greens-representation-bounded} decays exponentially fast as $M \to \infty$. As such, we may
replace $\Omega_M$ by its limit $\Omega$ in~\cref{eq:greens-representation-bounded}, yielding the following integral representation:
\begin{align}
  \label{eq:greens-representation}
  \tvarphi(\bx) = \mathcal{D}_{\Gamma}[\tvarphi](\bx) - \mathcal{S}_{\Gamma}[\nabla \tvarphi \cdot A \bn](\bx), \quad  \bx \in \Omega,
\end{align}
where $\Gamma := \partial \Omega = \Gamma_b \cup \Gamma_o \cup \Gamma_f$. 

It is worth mentioning that an integral representation on the unbounded curve~$\Gamma$, such as~\cref{eq:greens-representation}, can only be obtained for $\tvarphi$, and not for $\varphi$. In fact, in the absence of a complex-stretching, equation~\cref{eq:greens-representation} ceases to make sense because $\varphi$ does not decay, and the single-layer kernel exhibits a logarithmic growth, meaning that the integrals are not convergent (even in the conditional sense). This is a fundamental difference from the Helmholtz equation, where the decay of the Green's function and of the scattered field suffices to render integrals over unbounded interfaces conditionally convergent (but not absolutely so). 

Next, to obtain an integral equation for $\tvarphi$ over $\Gamma$, we evaluate the limit $\Omega \ni \bx \to \bx' \in \Gamma$ in~\cref{eq:greens-representation}. Taking into account the jump condition of the
the double-layer potential~\eqref{eq:DL-pot} across $\Gamma$~\cite{mclean2000strongly}, we arrive at %
\begin{align}
  \label{eq:greens-formula}
  \frac{\tvarphi(\bx)}{2} &=  D_{\Gamma}[\tvarphi](\bx) - S_{\Gamma}[\nabla \tvarphi \cdot A \bn](\bx), \quad \bx \in \Gamma,
\end{align}
where $S_{\Gamma}$ and $D_{\Gamma}$ are respectively the single- and double-layer operators,
defined as
\begin{subequations}
  \label{eq:calderon-operators}
  \begin{align}
  \label{eq:SL-operator}  
  S_{\Gamma}[\sigma](\bx) &:= \int_\Gamma \widetilde{G}(\bx,\by) \sigma(\by) \de s(\by) \quad\text{and}\\
  \label{eq:DL-operator}  
  D_{\Gamma}[\sigma](\bx) &:= {\rm p.v.}\! \int_\Gamma\left(\nabla_{\by} \widetilde{G}(\bx,\by) \cdot A(\by)\bn(\by) \right)
  \sigma(\by) \de s(\by), \quad \bx \in \Gamma.
\end{align}\label{eq:BIOS}\end{subequations}
The letters p.v. in front of the integral sign in~\cref{eq:DL-operator} refer to the fact that the integral must be interpreted in the sense of the Cauchy principal value.

Note that the double-layer kernel is given in explicit terms by
\begin{align}
  \label{eq:complex-double-layer-kernel}
\nabla_{\by}\widetilde{G}(\bx,\by) \cdot A(\by)\bn(\by) = \tau'(y_1)\frac{\partial G}{\partial\bn_{\by}}(\btau(\bx),\btau(\by)) = \frac{\tau'(y_1)(\btau(\bx)-\btau(\by)) \cdot \bn(\by)}{2\pi (\btau(\bx) - \btau(\by)) \cdot (\btau(\bx) - \btau(\by))},
\end{align}
which is simply a (complex) multiple of Laplace's free-space double-layer kernel evaluated at the complex points $\btau(\bx)$ and $\btau(\by)$. Furthermore, under the additional assumption that~$\Gamma$ and $\btau$ are smooth at $\bx$, it holds that
\begin{align}
\btau(\by) &= \btau(\bx) + J(\bx) (\by - \bx) + \mathcal{O}(|\bx - \by|^2), \\
\bn(\by) &= \bn(\bx) + \mathcal{O}(|\bx - \by|)
\end{align}
as $\by\to\bx$, where 
\begin{align}
    J(\bx) = \begin{bmatrix}
        \tau'(x_1) & 0\\
        0     & 1
    \end{bmatrix}
\end{align}
denotes the Jacobian matrix of $\btau$. Noting that $J(\bx) \bn(\bx) = \bn(\bx)$ for any $\bx \in \Gamma$, where we have used the fact that $J$ is either the identity (outside the PML) or $\bn$ points in the $x_2$ direction (inside the PML), we obtain that
\begin{align}
\nabla_{\by}\widetilde{G}(\bx,\by) \cdot A(\by)\bn(\by) = \frac{\tau'(x_1) (\bx - \by) \cdot \bn(\bx)}{J(\bx)(\bx - \by)\cdot J(\bx)(\bx - \by)} + \mathcal{O}(1) \quad \mbox{as} \quad \by \to \bx.
\end{align}
Thereofore, since $(\bx - \by) \cdot \bn(\bx) = \mathcal{O}(|\bx - \by|^2)$ as $\Gamma \ni \by \to \bx \in \Gamma$ for a smooth curve, we conclude that the double-layer kernel~\cref{eq:complex-double-layer-kernel} is a bounded function at $\by = \bx$, and hence the principal value in~\cref{eq:DL-operator} is not necessary in such cases. 

With the validity of the representation formula~\cref{eq:greens-representation}
established, and the availability of the free-space fundamental solution~\cref{eq:complex-Green-function}, we can
now derive a BIE for $\tvarphi$. Indeed, using the boundary conditions~\cref{eq:pml-free-surface-bc,eq:pml-solid-bc} to replace
$\nabla \tvarphi \cdot A \bn = \tau' \tvarphi \cdot \bn$ in favor of $\tvarphi$ in~\cref{eq:greens-formula}, we arrive at 
\begin{align}
  \label{eq:BIE-long}
  -\frac{\tvarphi(\bx)}{2} + D_{\Gamma}[\tvarphi](\bx) + \nu S_{\Gamma_f}\left[\tau'\tvarphi\right](\bx) &= S_{\Gamma_f}[f_1](\bx) + S_{\Gamma_b \cup \Gamma_o}[f_2](\bx), \quad \bx \in \Gamma.
\end{align}

To further simplify the notation, we introduce a ``global" Robin coefficient $\alpha : \Gamma \to \R$, as well as a ``global" source $f : \Gamma \to \R$ given by:
\begin{align}
    \alpha(\bx) = \nu\mathbbm{1}_{\Gamma_f}, \quad f(\bx) = f_1(\bx) \mathbbm{1}_{\Gamma_f} + f_2(\bx) \mathbbm{1}_{\Gamma_o \cup \Gamma_b},
\end{align}
where $\mathbbm{1}_{\Sigma}$ denotes the indicator function over a domain $\Sigma$. Equation~\cref{eq:BIE-long} is then rewritten in a more compact form as 
\begin{align}
  \label{eq:BIE}
  -\frac{\tvarphi(\bx)}{2} + D_{\Gamma}[\tvarphi](\bx) + S_\Gamma\left[\alpha\tau'\tvarphi\right](\bx) &= S_{\Gamma}[f](\bx), \quad \bx \in \Gamma.
\end{align}

In the subsequent sections, we will present the development of a Nystr\"om method for the numerical approximation of solutions to~\eqref{eq:BIE}.

\section{Numerical method}\label{sec:Nyst_meth}

In this section, we provide the details of the discretization scheme employed, as
well as the concrete choice of the (linear) PML function used along with the parameters involved. %
We also present an argument for the consistency of the truncated BIE which includes an integral equation for the error, that exhibits a right-hand-side that vanishes exponentially as $M\to\infty$. The numerical examples of~\Cref{sec:numerical-results} demonstrate that, as expected, the truncation errors decay exponentially with respect to the length of the PML layer.

In order to discretize~\cref{eq:BIE}, we proceed in two steps. First, after a choice of PML function $\tau$ in~\eqref{eq:tau-conditions} has been made, the infinite boundary $\Gamma$ is truncated to a finite curve $\Gamma_M = \{ \bx \in \Gamma : |x_1| < M \}$. Second, the resulting integral equation over $\Gamma_M$ is discretized using a Nyström scheme, wherein a composite $P$-point Gauss-Legendre quadrature rule is used to approximate the integrals over
$\Gamma_{M}$. A correction, which  accounts
for the singular and nearly-singular interactions arising in the numerical discretization of the underlying integral operators is obtained by means of a
Gauss-Kronrod adaptive quadrature; more details are given in Section~\ref{sec:discretization-error}. 

The particular choice of the BIE discretization here developed overcomes certain limitations of some existing Nystr\"om methods. For instance, the well-known spectrally accurate Martensen-Kussmaul approach~\cite[section 3.6]{colton1998inverse} becomes numerically unstable in the current setting due to the underlying kernel splitting  involving complex arguments~\cite{lu2014efficient}. Other approaches, such as the Density Interpolation Method~\cite{perez2019harmonic,faria2021general}, developed by some of the authors and collaborators, does not handle open curves such as the ones that arise in the current BIE formulation. The hybrid Gauss-trapezoidal quadrature rule of Alpert~\cite{alpert1999hybrid}, on the other hand, can in principle be used in the current setting and was the choice of spatial discretization used by Lu et al. in~\cite{lu2018perfectly}.

Throughout this section, we will use $\tvarphi_{M,h}$ to denote the approximate solution over $\Gamma_{M}$ using a $P$-point quadrature rule on elements of approximate size $h$, where the dependency on $P$ is omitted for notational simplicity. 

\subsection{Domain truncation}\label{sec:truncation}

To ensure a clear and concise presentation, we will focus exclusively on what we call ``linear PMLs'', wherein~$\tau$ is given by
\begin{align}
    \label{eq:linear-pml}
  \tau(x_1) &= x_1 + i 
  \begin{cases}
    c (x_1+a), & x_1 \leq -a, \\
    0, & |x_1| < a, \\
    c (x_1 - a) & x_1 \geq a,
\end{cases}
\end{align}
where $a>0$ is a parameter controlling the start of the complex stretching, while the parameter $c>0$ controls its slope. We know from~\cref{eq:modaldecomp} that $\tvarphi(x_1,x_2) = \varphi(\tau(\bx),x_2)$ admits the following
modal decomposition for $x_1 > L$:
\begin{align}
  \label{eq:modal-decomp-extension}
  \tvarphi(x_1,x_2) =  \displaystyle A_0^+u_0(x_2)\e^{-ck(x_1-a)} \e^{ik(x_1-L)} + \sum_{n\geq 1}A_n^+u_n(x_2)\e^{-\gamma_n(x_1-L)}\e^{-ic\gamma_n(x_1-a)}
\end{align}
with the coefficients $A_n^+$ given in~\cref{eq:modal-decomp-amplitude}. This shows that 
\begin{align}
    \tvarphi(x_1,x_2) \sim \mathcal{O}(\e^{-\min(ck,\gamma_1) x_1}) \quad \mbox{as} \quad x_1 \to \infty.
    \label{eq:decay-rate}
\end{align}
A similar argument holds for $x_1 \to -\infty$, and thus $\tvarphi$ decays exponentially as $|x_1| \to \infty$.

Given the exponential decay of $\tvarphi$, we replace the integral equation~\cref{eq:BIE} by its truncated version
\begin{align}
  \label{eq:BIE-truncated}
  -\frac{\tvarphi_{M}(\bx)}{2} + D_{\Gamma_{M}}[\tvarphi_M](\bx) + S_{\Gamma_{M}}\left[\alpha\tau'\tvarphi_M\right](\bx) &= S_{\Gamma_M}[f](\bx), \quad \bx \in \Gamma_M.
\end{align}
Subtracting~\cref{eq:BIE-truncated} from~\cref{eq:BIE}, it is straightforward to see that the truncation error $\mathcal{E}_M(\bx) := \tvarphi(\bx) - \tvarphi_M(\bx)$
satisfies the following integral equation:
\begin{align}
  \label{eq:BIE-truncated-error}
  -\frac{\mathcal{E}_M(\bx)}{2} + D_{\Gamma_{M}}[\mathcal{E}_M](\bx) + S_{\Gamma_{M}}\left[\alpha \tau' \mathcal{E}_M\right](\bx) &= -D_{\Gamma^c_{M}}[\tvarphi](\bx) - S_{\Gamma^c_{M}}\left[\alpha \tau' \tvarphi\right](\bx)
\end{align}
for $\bx \in \Gamma_M$, where $\Gamma_M^c = \Gamma \setminus \Gamma_M$. Because of the exponential decay of $\tvarphi$, it follows that the right-hand side is exponentially small in $M$ as $M \to \infty$, which shows the consistency of the truncation technique.  To further conclude that $\mathcal{E}_M$ is also exponentially small, however, a stability result of the type
\begin{align}
    || \mathcal{L}^{-1}_M ||_{L^2(\Gamma_M) \to L^2(\Gamma_M)} \lesssim M^q, \quad    \mathcal{L}[\sigma] := - \frac{\sigma}{2} + D_{\Gamma_M}[\sigma] + S_{\Gamma_{M}}[\alpha \tau' \sigma]
\end{align}
is needed as $M \to \infty$, with $q > 0$ independent of $M$. Although such a result is lacking in our analysis, the numerical tests of~\Cref{sec:numerical-results} suggest that the resulting truncation errors do indeed decay exponentially fast as $M$ increases. 

\begin{remark}[PML parameters]
While it may be tempting to take the parameter $c$ large in~\cref{eq:linear-pml} so as to improve the exponential decay of $\tvarphi$, consideration of discretization errors ($\tvarphi$ has to be approximated on a grid) will typically put bounds on the size of $c$. Even if we assume, as is the case in this paper, that the same grid size is used inside and outside the PML layer, the question of determining the ``optimal'' PML parameters is subtle and depends closely on the discretization scheme employed, for one seeks to balance truncation and discretization errors. In this paper, we simply take $c=1$, and avoid questions related to an optimal choice of PML parameters. 
\end{remark}

\subsection{Singular and nearly-singular integration}
\label{sec:discretization-error}

With a truncated curve $\Gamma_M$ at hand, we are now in position to discretize the integral operators in~\cref{eq:BIE-truncated}. To this end, we assume that the curve $\Gamma_M$ is covered by $N$ non-overlapping possibly curved elements $\left\{ \mu_n \right\}_{n=1}^N$ of approximate size $h$; i.e. 
\begin{align}
  \Gamma_M = \bigcup_{n=1}^N \overline{\mu}_n.
\end{align}
For each element $\mu_n$, a $\mathcal{C}^\infty$ bijection $\chi_n : [0,1] \to \overline{\mu_n}$ is assumed known.

On each $\mu_n$, a $P$-point Gauss-Legendre quadrature rule $\mathcal{Q}_n = \left\{ \left( \bx_{p,m}, w_{p,m} \right) \right\}_{p=1}^P$ is constructed by mapping the quadrature on the reference segment $[0,1]$ to $\overline{\mu_m}$ using $\chi_m$. A composite quadrature $\mathcal{Q}$ for $\Gamma_{M}$ is then obtained by collecting the element quadratures:
\begin{align}
  \mathcal{Q} = \left\{ \mathcal{Q}_n \right\}_{n=1}^N = \left\{ \left\{ \left( \bx_{p,n}, w_{p,n} \right) \right\}_{p=1}^P \right\}_{n=1}^N.
\end{align}

Given an element index $p$ and a local quadrature index $n$, it will be sometimes convenient to use a linear index $i = n(P-1) + p$, $1 \leq i \leq NP$, to refer to the node $\bx_{p,n}$ as $\bx_i$. With this linear index, the boundary quadrature can be alternatively written as 
 $ \mathcal{Q} = \left\{ \left( \bx_{i}, w_{i} \right) \right\}_{i=1}^{NP}.$

Using a Nyström discretization, the integral operators in~\cref{eq:BIE} are approximated using the quadrature $\mathcal{Q}$, and equality is enforced (strongly) at the quadrature nodes. This yields the following discrete linear system for $\tvarphi_{i} \approx\tvarphi_M(\bx_i)$:
\begin{align}
    \label{eq:discrete-BIE}
    -\frac{\tvarphi_i}{2} + \sum_{j=1}^{NP} K_{i,j}\tvarphi_j  =
    \sum_{j=1}^{NP} G(\bx_i,\bx_j) f(\bx_j)w_j,    \quad  1 \leq i \leq NP,
\end{align}
where  $K_{i,j} := K(\bx_i,\bx_j)w_j$ with
\begin{equation}\label{eq:kernel_disc}
K(\bx,\by) :=  \nabla_{\by} \widetilde{G}(\bx,\by) \cdot A(\by)\bn(\by) + \alpha(\by)\tau'(\by)\widetilde{G}(\bx,\by).    
\end{equation}

It is well-known that, for a given element $\mu_n\subset\Gamma_M$, the quadrature rule $\mathcal{Q}_n$ need to be corrected for points $\bx_i$ close to (or on) $\mu_n$ due to the singular nature of the kernels in the integral operators. In other words, the discrete problem~\cref{eq:discrete-BIE}, needs to be modified to properly account for the singular and nearly-singular interactions, where the Gauss-Legendre quadrature becomes inaccurate due the lack of smoothness of the integrand. In this paper, we opt for a simple correction method whereas $\tvarphi$ is replaced by its Lagrange interpolant on the element $\mu_n$, and the integration against the Lagrange interpolant is performed using an adaptive Gauss-Kronrod quadrature. In more detail, for a target point $\bx_i$ ``close'' to the element $\mu_n$, where close is measured relative to the local mesh size, the matrix entries $K_{i,j_p}$ obtained from the regular quadrature are replaced by the following integrals
\begin{align}
    \label{eq:singular-weights}
    \widetilde{K}_{i,j_p} = \int_{\mu_n} K(\bx_i,\by) \ell_{n,p}(\by) \de s_{\by},\quad 1\leq p\leq P,
\end{align}
 where $j_p=n(P-1)+p$ and $\ell_{m,p}$ is defined as $\ell_{n,p}(\chi_n(s)) = \hat{\ell}_p(s)$,   $0 \leq s \leq 1$ with $\hat{\ell}_p$ being the $p$-th Lagrange basis over the reference $P$-point Gauss-Legendre quadrature $\left\{ s_p, w_p \right\}_{p=1}^P$; i.e., $\hat{\ell}$ is the $P-1$ degree polynomial satisfying $\hat{\ell}_p(s_q) = \delta_{pq}$ for $1 \leq q \leq P$. In practice, the integrals in~\cref{eq:singular-weights} are computed using an adaptive quadrature rule to a prescribed absolute tolerance; in all the examples presented in this paper, we set the integration tolerance to $10^{-10}$. The same correction approach is taken to accurately evaluate the right-hand-side of~\eqref{eq:discrete-BIE} which involves the single-layer operator.

The errors introduced by the corrected quadrature approach for the singular and nearly-singular entries can be estimated as follows. First, using the definition of $\widetilde{K}_{i,j_p}$ in~\cref{eq:singular-weights} and assuming that $\widetilde\varphi(\bx_{j_p})=\widetilde\varphi_{j_p}$, the quadrature error for the element $\mu_n$ with a target point $\bx_i$ can be written as
\begin{align}
    \label{eq:singular-quadrature-split}
    \int_{\mu_n} K(\bx_i,\by) \tvarphi(\by) \de s_{\by} - \sum_{p=1}^P \tvarphi_{j_p} \widetilde{K}_{i,j_p} =
    \int_{\mu_n} K(\bx_i,\by) \left( \tvarphi(\by) - \sum_{p=1}^P \tvarphi_{j_p} \ell_{n,p}(\by) \right) \de s_{\by}.
\end{align}
From standard interpolation results (e.g.~\cite[Theorem 1]{strang1972approximation}), we have that 
\begin{align}
    \max_{\by \in \tau_m} \left|  \tvarphi(\by) - \sum_{j=1}^P \tvarphi_{j_p} \ell_{n,j}(\by)  \right| < C h^P,
\end{align}
for some constant $C$ which depends on the derivatives up to order $P$ of $\tvarphi$ and the element parametrization $\chi_n$. 
Provided $K(\bx_i,\cdot)$ is absolutely integrable over $\mu_n$, as in all examples of this paper, we may then write the following bound:
\begin{align}
    \label{eq:lagrange-integration-error}
    \left| \int_{\mu_n} K(\bx_i,\by) \left( \tvarphi(\by) - \sum_{p=1}^P \tvarphi_{j_p} \ell_{m,j}(\by) \right) \de s_{\by} \right| \leq C h^{P} \int_{\mu_n} | K(\bx_i,\by) | \de s_{\by}.
\end{align}

It is now important to distinguish between the single and double-layer kernels making up $K$, as
their singularities as $\by\to\bx$ differ. The single-layer kernel has a singularity of
logarithmic type; this means, in particular,
that its contribution to the integral on the right-hand side of~\cref{eq:lagrange-integration-error}
can be bounded by $h |\log(h)|$. The behaviour of the double-layer kernel is slightly more subtle. If
$\bx_i \in \mu_m$, the double-layer kernel is in fact bounded, and its contribution to the integral on the right-hand side of~\cref{eq:lagrange-integration-error}
can be bounded by $h$, yielding a logarithmic improvement over the single-layer
error. If $\Gamma \ni \bx_i \not \in \mu_m$ with $\mbox{dist}(\bx_i,\mu_m) \sim h$, however, two cases arise depending on the regularity of the curve $\Gamma_M$ (which is assumed to be at least Lipschitz). These are
\begin{align}
  \int_{\tau_m} \frac{(\btau(\bx_i) - \btau(\by))\cdot \bn(\by)}{|\bx_i - \by|^2}  \de s_{\by} \sim 
  \begin{cases}
    \mathcal{O}(h) \quad \mbox{if $\Gamma$ is $\mathcal{C}^2$}\\
    \mathcal{O}(1) \quad \mbox{otherwise}
  \end{cases}
\end{align}
Since $K$ in~\cref{eq:kernel_disc} involves both the single- and
double-layer kernels, we expect the error in~\eqref{eq:lagrange-integration-error} to be proportional to $h^{P+1} |\log(h)|$ for smooth curves, or $h^P$ for curves featuring corners. The numerical experiments presented in next section, closely match these estimates.

\section{Numerical results}\label{sec:numerical-results}

In this section, we present a variety of numerical results designed to showcase the capabilities of the proposed methodology,
as well as to test the effect of the  parameters involved in the domain truncation and discretization. All linear systems are solved iteratively by means of GMRES~\cite{saad1986gmres}, with a relative residual tolerance of $10^{-10}$. 

In the examples that follow, the reported errors are measured at a finite set of test points $\boldsymbol{X} $ located either on the boundary $\Gamma$ (in which case they coincide with the quadrature nodes) or in the volume $\Omega$. In detail, errors are measured via 
\begin{align}
  E_{M,h} = \max_{\bx \in \boldsymbol{X}} |\varphi_{M,h}(\bx) - \varphi_{\rm ref}(\bx) |,
\end{align}
where $\varphi_{\rm ref}$ denotes a reference solution which is the exact PDE solution when it is available or a
 solution obtained using a fine discretization of $\Gamma_M$ with a large truncation parameter $M$. When $\bx=\bx_j\in\boldsymbol{X}\cap\Gamma_M$ is a quadrature node, the value $\varphi_{M,h}(\bx)$ is taken from  $\widetilde\varphi_{j}$ which is in turn obtained by solving  the underlying linear system. Otherwise,  when $\bx \in\boldsymbol{X}\cap \Omega$, we use the domain quadrature $\mathcal{Q}$ to evaluate an approximation to~\cref{eq:greens-representation}:
\begin{align}
    \varphi_{M,h}(\bx) = \sum_{j=1}^{NP} K(\bx,\bx_j) \tvarphi_j -\widetilde{G}(\bx,\bx_j)f(\bx_j)w_j. 
\end{align}

\subsection{Wave-maker}\label{subsec:wavemaker}

The first example problem that we consider is that of a wavemaker whose geometry consists of the semi-infinite strip $[0,\infty) \times [-1,0]$. For validation purposes, we construct an exact solution by summing the right-going propagating mode and the first evanescent mode in~\cref{eq:modaldecomp}:
\begin{align}
    \label{eq:ref-solution-wavemaker}
  \varphi_{\rm ref}(x_1,x_2) = \frac{\cosh(k(x_2+1))}{\cosh(k)} \e^{i k x_1} + \frac{\cos(\gamma_1(x_2+1))}{\cos(\gamma_1)} \e^{-\gamma_1 x_1},
\end{align}
where $k$ and $\gamma_1$ are implicitly defined by~\cref{eq:dispersion-relation} and~\eqref{eq:disp-evanescent}, respectively. 

The Neumann trace of $\varphi_{\rm ref}$ is given as a source on the left boundary $\{ 0 \} \times[-1,0]$, and a numerical approximation $\tvarphi_{M,h}$ is found by solving the linear system~\cref{eq:discrete-BIE}. For this example, the error is measured at all the quadrature nodes of $\Gamma_M$; for points inside the PML, we use $\tvarphi_{\rm ref}(x_1,x_2) = \varphi_{\rm ref}(\tau(x_1,x_2),x_2)$ as the exact solution. The PML function is given by~\cref{eq:linear-pml}, with the PML strength $c=1$. The value of $a$, denoting the starting point of the PML, is indicated in each example. 

\begin{figure}[ht!]
  \centering
  \vspace{-5pt}
  \begin{subfigure}[t]{0.49\linewidth}
    \includegraphics[width=1\textwidth]{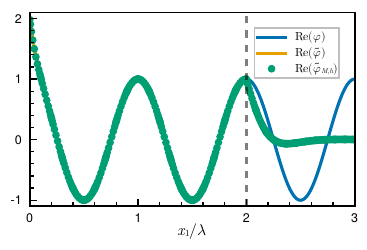}
    \caption{Exact solution $\varphi_{\text{ref}}$, its analytic extension
    $\tvarphi_{\text{ref}}$, and the numerical approximation $\tvarphi_{\ell,h}$. The vertical dashed line indicates the start of the PML layer.}
    \label{fig:wavemaker-comparison-exact-numerical}
  \end{subfigure}\hfill
  \begin{subfigure}[t]{0.49\linewidth}
    \includegraphics[width=1\textwidth]{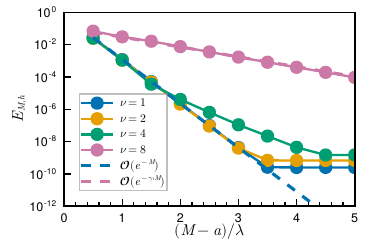}
    \caption{Convergence to exact solution as the truncation length $M$ is increased for different values of the parameter $\nu$. The dashed lines are the expected reference slopes.}
    \label{fig:convergence-wavemaker}
  \end{subfigure}
  \caption{Solution to the wave-maker problem (left) and convergence for various values of $\nu$ (right).}
  \vspace{-10pt}
\end{figure}

\Cref{fig:wavemaker-comparison-exact-numerical}  displays the solution of the aforementioned problem on the free surface for the $\nu= d \omega^2/g  = 1$, which yields $k \approx 1.2$, and $\gamma_1 \approx 2.8$. The PML layer starts at $a=2\lambda$ (indicated by a dashed vertical line), where $\lambda = 2\pi/k$ is the surface wavelength. Since we are interested in assessing the effect of the truncation on the error in this example, we fix $h = \lambda/10$ and $P=10$; i.e., the resulting domain discretization consists of $100$ points-per-wavelength, which in this example suffices to achieve 10 digits of accuracy in the approximate BIE solution. A visual comparison between the exact solution, its analytical extension, and the numerical approximation $\tvarphi_{M,h}$ is shown in~\cref{fig:wavemaker-comparison-exact-numerical}, where it is apparent that $\tvarphi_{M,h}$ converges to $\tvarphi$ on all of $\Gamma_M$, and therefore to $\varphi$ for $x_1 < a = 2\lambda$ (we only show here the values at $\Gamma_f$). For~\cref{fig:wavemaker-comparison-exact-numerical}, the domain is truncated at $M=3\lambda$. 
To study the rate of convergence as $M$ increases, we then present \Cref{fig:convergence-wavemaker} which displays the maximum error $E_{M,h}$, measured again on all the quadrature nodes, as a function of the truncation parameter $M$ for various values of $\nu$. (As in the previous example, we use a fine domain discretization to suppress the effect of BIE discretization errors.) We vary $M$ from $a$ to $a + 5 \lambda$ and display $E_{M,h}$ as a function of the ``normalized'' PML length $\ell = (M-a)/\lambda$. Given the exponential decay rate of $\varphi$ in~\cref{eq:decay-rate}, we expect a truncation error of $\mathcal{O}(\e^{-\min(k,\gamma_1)M})$ which can be more conveniently expressed in terms of the dimensionless PML length as $\mathcal{O}(\e^{-\min(2\pi,\gamma_1 \lambda) \ell} )$.  The expected decay rate is clearly observed in~\cref{fig:convergence-wavemaker}, where for small $\nu$ the error is dominated by the propagating mode (and thus exhibits a $\mathcal O(\e^{-kM})$ rate), while for large values of $\nu$ we observe a slower exponential decay due to the evanescent mode (which decays as $\mathcal{O}(\e^{-\gamma_1 M})$). The plateau observed around $10^{-10}$ is an indication of errors stemming from the BIE solution (which include, for instance, the iterative solver residual and the quadrature rule errors) becoming dominant.

\begin{figure}[ht!]
  \centering
  \includegraphics[width=0.5\textwidth]{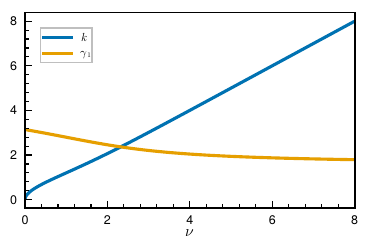}
  \caption{Wavenumber $k$ versus $\gamma_1$ as a
  function of $\nu$.}
  \vspace{-10pt}
  \label{fig:propagative-vs-evanescent}
\end{figure}

As observed in~\cref{fig:convergence-wavemaker}, the exponential decay rate of truncation errors appears to deteriorate as $\nu$ increases. This is easy to understand if we consider~$k$ and~$\gamma_1$ as functions of $\nu$; both functions are plotted in~\cref{fig:propagative-vs-evanescent}. As can be seen (and it easily follows from~\cref{eq:dispersion-relation,eq:disp-evanescent}) we have that $k \to \nu$, $\gamma_1 \to \pi/2$ in the limit $\nu \to \infty$. As such, the truncation of the evanescent mode becomes dominant for sufficiently large~$\nu$, and a large PML layer may be needed to obtain a small truncation error. While not a major drawback of the method for reasonably small $\nu$ values --- the method still converges exponentially fast for any fixed $\nu$ --- we propose next a possible modification to the PML function~\cref{eq:linear-pml} to improve the truncation errors in the presence of slowly-decaying evanescent modes. 

\begin{figure}[ht]
  \centering
  \begin{subfigure}{0.49\linewidth}
    \includegraphics[width=\textwidth]{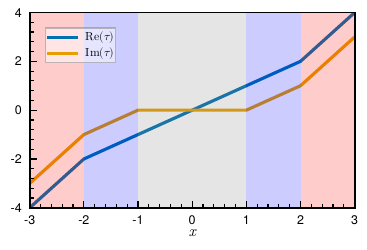}
    \caption{Two-layered PML with $a=1$, $b=2$, $c=1$, $\nu=2$}
    \label{fig:two-layer-pml}
  \end{subfigure}\hfill
  \begin{subfigure}{0.49\linewidth}
    \includegraphics[width=1\textwidth]{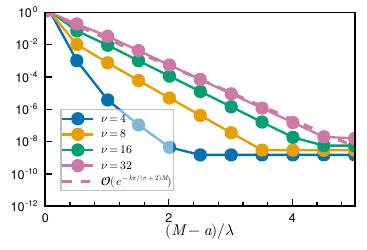}
    \caption{Convergence for large values of $\nu$ with two-layer PML.}
    \label{fig:convergence-stretching-vary-depth}
  \end{subfigure}
  \label{fig:convergence-modal-solution}
  \caption{Two-layered PML function (left) and its convergence properties (right)}
\end{figure}

Inspired by~\cite{Vial:12}, we improve our complex stretching to account not only for the propagating surface-wave mode but also for the evanescent modes, hence ensuring the effectiveness of the PML for a larger range of $\nu$ values. The main idea is to employ two different absorbing layers; a complex scaled layer, as before, to absorb the propagating mode, followed by a real stretching to accelerate the convergence of the evanescent modes. In  detail, we replace the PML function~\cref{eq:linear-pml} by:
\begin{align}
    \label{eq:two-layer-pml}
  \tau_{\text{two-layer}}(x_1) &=  
  \begin{cases}
    \tau(x_1), & |x_1| \leq b, \\
    \tau(\nu x_1) - \tau(\nu b) + \tau(b) & |x_1| > b, \\
\end{cases}
\end{align}
where $b > a$ is a new parameter controlling the start of the real-stretching layer. The two-layered PML is depicted in~\cref{fig:two-layer-pml}. 
It is important to note that the real stretching increases the frequency of the oscillatory mode, and therefore should only be applied \emph{after} the complex stretching layer.

Given the start of the PML layer $x_1=a$, and the truncation parameter $M$, we chose $b$ by the following simple rationale. Supposing the propagating mode decays by~$\e^{-ck(b-a)}$, while the slowest decaying evanescent mode evanescent mode decays by $\e^{-\gamma_1 \nu (M-b)}$, requiring both errors to be equal yields
\begin{align}
 b = \frac{cka + \gamma_1 \nu M}{ck+\gamma_1 \nu}.
 \label{eq:b-formula}
\end{align} 
\Cref{eq:b-formula} provides a way of choosing $b$ given the other problem parameters. Inserting this expression into the decay rate of the propagating mode, $\e^{-ck(b-a)}$, and expressing the result in terms of the truncation parameter $M$, we obtain:
\begin{equation}
    \e^{-ck(b-a)} \sim \mathcal{O}(\e^{-ck \mu M}), \quad \mbox{where} \quad \mu = \frac{\gamma_1 \nu}{ck + \gamma_1 \nu}.
    \label{eq:two-layer-decay}
\end{equation}

\Cref{eq:two-layer-decay} provides the expected decay rate of $\tvarphi$ when the two-layer PML function~\cref{eq:two-layer-pml} is used. The constant $\mu < 1$ in~\cref{eq:two-layer-decay} provides an idea of the deterioration compared to~\cref{eq:decay-rate} for the cases where $ck < \gamma_1$ (i.e., cases for which the decay is dominated by the propagating mode). Note that, as $\nu \to \infty$, we have $\mu \to \pi / (2c + \pi)$, meaning the decay rate is at worst a constant factor slower than~$\e^{-ck\mu}$; for $c = 1$, for example, we have $\mu \approx 0.61$ as $\nu \to \infty$. As can be seen in~\cref{fig:convergence-stretching-vary-depth}, the proposed idea significantly improves the truncation error, even for values of $\nu$ as large as $32$. Furthermore, we observe that as $\nu$ increases the decay appears to match the expected $\mathcal{O}(\e^{-k \pi/(\pi+2)M})$ asymptote, shown as a dashed line in~\cref{fig:convergence-stretching-vary-depth}. 

\begin{figure}
  \centering
  \includegraphics[width=0.5\textwidth]{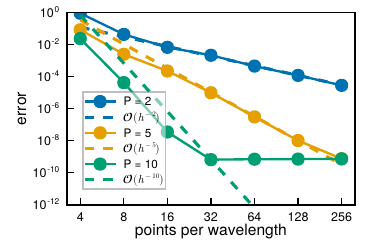}
  \caption{Convergence to the modal solution for $\nu=2$, and $P=2,5,10$, as the number of
  points per wavelength is increased. The expected convergence rate of order $P$
  described in Section~\ref{sec:discretization-error}, shown as dashed lines of the same color, is achieved.}
  \label{fig:mesh-convergence}
  \vspace{-15pt}
\end{figure}

Finally, although not the main focus of this paper, we present in~\cref{fig:mesh-convergence} convergence results with respect to mesh refinement. 
Following the discussion in Section~\ref{sec:discretization-error}, we expect a convergence of order $P$ for the wavemaker, which is what is observed~\cref{fig:mesh-convergence}. (For this example we took a PML of length $4 \lambda$, where $\lambda$ is the wavelength, which was enough to ensure sufficiently small truncation errors.)

\subsection{Scattering problem}

Having in the previous sections validated the proposed algorithm and its implementation, we here focus on the classical scattering problem. In the  examples that follow we consider an incident field, corresponding to a right-going propagating mode, given by 
\begin{align}
  \varphi^{\rm inc}(x_1,x_2) = \frac{\cosh(k(x_2+1))}{\cosh(k)} \e^{i k x_1}.
  \label{eq:inc-field}
\end{align}
Expressing the total field as
\begin{align}
  \varphi^{\rm tot} = \varphi^{\rm inc} + \varphi,
\end{align}
we seek to find the scattered field $\varphi$ satisfying the radiation condition~\cref{eq:radiation-condition}. It follows from the problem linearity that $\varphi$
solves~\cref{eq:water-waves-system} with $f_1 = 0$ and $f_2 = -\nabla \varphi^{\rm inc}
\cdot \bn$; i.e.~$\varphi$ satisfies the following boundary condition:
\begin{align}
  \varphi = - \nabla \varphi^{\rm inc} \cdot \bn, \quad \bx \in \Gamma_o \cup \Gamma_b.
\end{align}

We next consider three different settings for the scattering problem consisting of (i) a finite number of fully submerged jellyfish-like obstacles over a variable topography, (ii) a constant-depth topography with two piercing obstacles, and (iii) a topography containing a step profile. In all of the examples considered in this section we use a mesh of size $\lambda/30$ and the parameter values $P=5$, and $\nu=4$ (which yields a wavelength $\lambda \approx 1.57$); the PML parameters are fixed as $a=3\lambda$, $c=1$, and the domain is truncated at $M=4\lambda$.

\subsubsection{Submerged obstacles}\label{sec:submerged-obs}

In our first scattering problem, we consider a bottom topography given by the following height function:
\begin{align}
    \mathrm{topography}(s) = (s,-1 + \e^{-s^2/4} \sin(2\pi s /\lambda)),\quad s\in\R,
\end{align}
which gives rise to  a ``bump-like'' perturbation to a constant depth
\footnote{Although this is not a compact perturbation of a constant depth, the perturbations decay exponentially fast, and are numerically indistinguishable from a compact perturbation while being easier to construct.}
. In addition to the bump-like topography, we consider three jellyfish-like obstacles whose boundaries are given by the following parametrization (after appropriate rotation, translation, and scaling):
\begin{align}
    \mathrm{jellyfish}(s) = r(s)(\cos(s),\sin(s)),\quad s\in\R,
\end{align}
where $r(s) = 1+0.3\cos(4s+2\sin(s))$.

\begin{figure}[ht]
  \centering
  \vspace{-30pt}
  \includegraphics[width=1\textwidth]{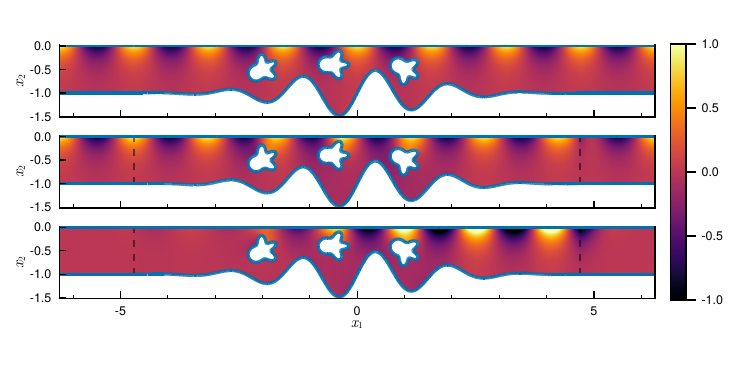}
  \vspace{-30pt}
  \caption{Incident (top), total (middle), and scattered (bottom) fields for three jellyfish-like obstacles over a variable topography. The dashed lines indicate the start of the PML layer, after which we observe a decay of the scattered field (bottom) due to the complex stretching.}
  \label{fig:jellyfish-fields}
  \vspace{-10pt}
\end{figure}

The resulting geometry is displayed in~\cref{fig:jellyfish-fields} along with the incident, total, and scattered, fields. The damping of the scattered field inside the PML-layer is clearly visible in~\cref{fig:jellyfish-fields} (bottom), especially for $x_1 > 2\lambda$ (i.e., passed the dashed lines). The main effect of the topography and obstacles is to cause a phase shift in the transmitted wave, as can be observed by comparing the incident field (top) to the total field (center) in~\cref{fig:jellyfish-fields}.

Because no exact solution is known for this scattering problem, we verified convergence of our Nystr\"om method with respect to the mesh size by performing a self-convergence test, where we evaluate $\varphi$ at a uniform mesh of the rectangle $[1.5\lambda,2.5\lambda] \times [-0.75 \times -0.25] \subset \Omega$, and compare the results to a reference solution obtained using $256$ points-per-wavelength BIE domain discretization. We expect from the discussion in Section~\ref{sec:discretization-error} the error to decrease as $\mathcal O(h^{6})$ as $h\to0$ (this is a smooth geometry), which appears to agree with the observed convergence rate in~\cref{fig:convergence-jellyfish}. 

\subsubsection{Piercing obstacles}\label{sec:piercing-geo}
While we have insofar considered only fully submerged obstacles (i.e., the free surface coincides with $\R \times \{0\}$), piercing objects can be handled in a similar fashion. From an application viewpoint, piercing obstacles are particularly interesting since many relevant scatterers (ships, docks, buoys) lie precisely on the free surface. For simplicity, we take here a geometry given by two piercing semi-circles of radius $r=1/4$ centered at $x_1 = \pm x_c$, with $x_c = r + 1/2$. This geometry creates a free surface with three disconnected components: $(-\infty, -1] \times \left\{0 \right\}$, $[-0.5,0.5] \times \left\{0 \right\}$, and $[0.5, \infty) \times \left\{0 \right\}$. As we will see later in Section~\ref{sec:resonant-freq}, the part of the surface disconnected from infinity can act as a resonator for certain values of the Robin parameter~$\nu$.

\begin{figure}[ht!]
  \centering
  \vspace{-10pt}
  \includegraphics[width=1\textwidth]{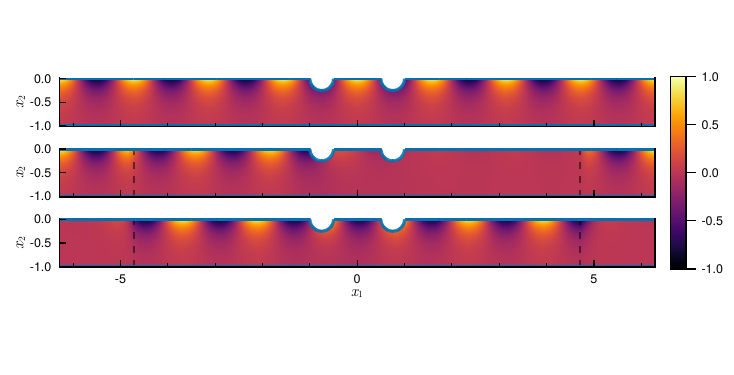}
  \vspace{-40pt}
  \caption{Incident (top), total (middle), and scattered (bottom) fields in the presence of two piercing bodies. The dashed vertical lines indicate the start of the PML layer.}
  \label{fig:piercing-fields}
\end{figure}

As for the other geometries considered above, we display in~\cref{fig:piercing-fields} the incident, total, and scattered fields. Unlike the submerged obstacle, we observe a large scattered field both on the left and on the right of the obstacles. Again, as expected, the scattered field decays within the PML (bottom figure), both on the right and on the left of the obstacles. 

\begin{figure}[ht]
\begin{subfigure}[t]{0.49\textwidth}
  \includegraphics[width=1\textwidth]{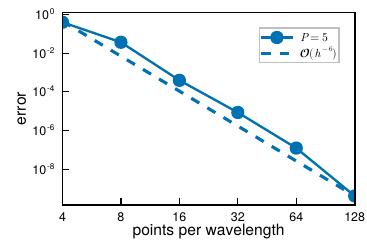}
  \caption{Smooth geometry (jellyfish)}
  \label{fig:convergence-jellyfish}
\end{subfigure}
\hfill
\begin{subfigure}[t]{0.49\textwidth}
    \includegraphics[width=1\textwidth]{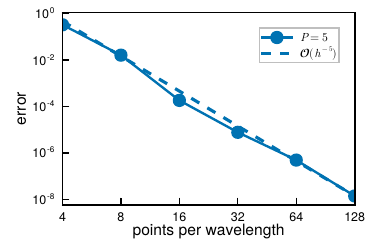}   
      \caption{Non-smooth geometry (piercing obstacle)}
    \label{fig:piercing-convergence}
\end{subfigure}

\caption{Self-convergence of scattering problem for the geometries shown in~\cref{fig:jellyfish-fields}(a) and~\cref{fig:piercing-fields}(b). A reference solution with $256$ points-per-wavelength was used when computing the error.}
\end{figure}

Finally, we perform the same self-convergence test described in Section~\ref{sec:submerged-obs}, where we evaluate the maximum error of the solution on a test domain inside of $\Omega$. The convergence with respect to the mesh refinement is shown in~\cref{fig:piercing-convergence}, where the reference line has the slope corresponding to $\mathcal O(h^{-5})$. Since the piercing geometry features corners, we expect such behavior according to Section~\ref{sec:discretization-error}, although a detailed analysis of the discretization errors in such cases is likely to be significantly more challenging due to the abrupt switch in the boundary condition arising at the junction points connecting the piercing obstacle and the free surface. 

\subsubsection{Step topography}

The final scattering geometry considered is that of a topography with a different depth as $x_1 \to \pm \infty$. In particular, we use here the following smooth step-like topography function:
\begin{align}\label{eq:step_top}
    \textrm{step}(s) = \left(s,-1 + \frac{\Delta d}{2} (1+\tanh(5s))\right), \quad s\in\R,
\end{align}
where, for this particular example, we take~$\Delta d = 0.75$ so that the shallow region on the right-hand side of the bottom has a depth of~$0.25$. It is worth mentioning that this is a challenging problem for methods based on the use of the standard water waves Green's function because the Neumann trace of the Green's function fails to vanish on an unbounded curve, making it necessary to introduce suitable approximations to achieve a BIE that can be discretized and numerically solved (see e.g.~\cite{porter2000water}).

This type of topography does not enter the framework introduced in~\Cref{sec:prob_setup} as it does not give rise to a computational domain that is a compact perturbation of a unit depth waveguide. Incidentally, the source term~$f_2$ is not compactly supported, since the Neumann trace of the incident right-going wave does not vanish at $x_2 = -0.25$ for large values of $x_1$ (i.e., it does not satisfy the boundary condition on the right-hand side of the bottom boundary). Upon closer inspection, however, it is straightforward to see that both $\tvarphi^{\rm inc}(\bx) = \varphi^{\rm inc}(\btau(\bx))$ and  its Neumann trace ${\bn} \cdot \nabla\tvarphi^{\rm inc}$, decay exponentially fast as $x_1 \to \infty$ (recall the definition of $\varphi^{\rm inc}$ in~\cref{eq:inc-field}), implying that the analytic extension of~$f_2$, although not compactly supported, decays exponentially fast as $x_1 \to \infty$. Although only a formal argument, the exponential decay of the source term means that we can expect our numerical method to behave in this case similarly to the compactly supported case.  

\begin{figure}[ht!]
  \centering
  \vspace{-10pt}
  \includegraphics[width=1\textwidth]{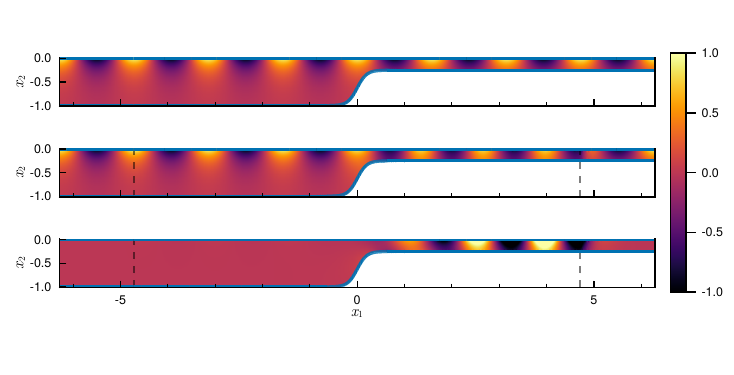}
  \vspace{-40pt}
  \caption{Incident (top), total (middle) and scattered (bottom) field for step-like bottom topography given by the curve parametrization in~\eqref{eq:step_top}.}
  \vspace{-15pt}
  \label{fig:step-fields}
\end{figure}

As it turns out, this is indeed what we observe in~\cref{fig:step-fields}, where we display the incident, total, and scattered field for the step-like geometry. The decay of the scattered field $\tvarphi$ inside the PML layer becomes apparent (bottom image), and we see that for this value of the Robin parameter $\nu$ ($\nu=4$ in this case), the scattered field is mostly concentrated on the shallow region. The change in the wavelength in the total field (middle image) is clearly visible on the right of the step, where the effective frequency is increased in the shallower region. While further theoretical and numerical investigations are needed to fully validate the effectiveness of the proposed technique for step-like topographies, we believe the results obtained are promising, in particular given the simplicity of the method. 

\subsection{Resonant frequencies}\label{sec:resonant-freq}

In this final set of examples, we consider the problem of finding resonant frequencies to the water-waves problem. More precisely, we seek to find the complex values of the Robin parameter $\nu$ such that~\cref{eq:water-waves-system} has a nontrivial solution for $f_1=f_2=0$. Using the boundary integral formulation~\cref{eq:BIE-long}, this corresponds to finding eigenvalues $\nu$ and eigenfunctions $\tvarphi_\nu$ solving the following (generalized) linear eigenvalue problem:
\begin{align}
  \label{eq:BIE-eigenvalue}
  -\frac{\tvarphi_\nu(\bx)}{2} + D_{\Gamma}[\tvarphi_\nu](\bx) =  -\nu S_{\Gamma_f}\left[\tau'\tvarphi_\nu\right](\bx), \quad \bx \in \Gamma.
\end{align}

Note that we do not allow for $\tau$ to depend on $\nu$, which would make for a nonlinear eigenvalue problem. In what follows, we fix the PML  parameters to be $a=2.5$, $c=1$, and take $M = 5$; for the discretization parameters, we set $P=5$ and take a mesh size of $h=0.2$. 

\begin{figure}[ht!]
  \centering
  \includegraphics[width=1\textwidth]{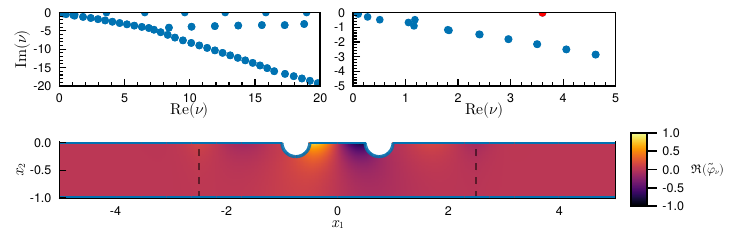}
  \vspace{-20pt}
  \caption{Complex resonances of a piercing geometry with two semi-circles. The top-left image shows the spectrum on a large portion of the complex plane, while the top-right image provides a zoomed-in version where we can clearly see the presence of an ``almost real'' eigenvalue, depicted in red. The bottom image displays the eigenfunction corresponding to the red eigenvalue in the top-right image. }
  \vspace{-10pt}
  \label{fig:double-piercing-spectrum}
\end{figure}

\Cref{fig:double-piercing-spectrum} displays the spectrum corresponding to the piercing geometry of Section~\ref{sec:piercing-geo}, as well as the eigenfunction corresponding to an eigenvalue with a small imaginary part. The eigenvalues on a large region of the complex plane enclosing $\{z\in\C: 0\leq \real z\leq 20, -20\leq\imag z\leq 0\}$ are shown on the top-left figure. There we can clearly observe that a part of the spectrum follows a straight line of approximate slope $-1$; these eigenvalues are typically PML-dependent and correspond to the complex stretching of the continuous spectrum $\R^+$ in the original problem. The bifurcation at $\real \nu \approx 10$ has also been observed in~\cite[e.g. Figures 5,6]{treyssede2014finite}. The more interesting and physically relevant eigenvalues  appear at nearly regular intervals close to the real axis. These correspond to eigenfunctions that are localized near the obstacles, and are not noticeable affected by variations in the PML parameters. We show in the top-right image a zoomed-in version of the spectrum, and mark in red one of the aforementioned eigenvalues with a small imaginary part. Finally, the bottom image displays the (normalized) eigenfunction corresponding to the highlighted red eigenvalue. We see that portion of the free surface between the two-obstacles acts as a trapping region, where localized eigenfunctions can exist. 

\section{Open theoretical problems}\label{sec:open-problems}

This work has opened up several challenging theoretical questions that were not addressed in this paper. In what follows we mention a few of them.

Firstly, the logarithmic growth of complex-scaled Green function~\cref{eq:complex-Green-function} makes it difficult to pose a functional framework for analyzing the well-posedness of~\cref{eq:BIE} in the usual Sobolev spaces. A possible way to circumvent such theoretical difficulties is to consider instead a Green's function $\widetilde{G}^D$ satisfying a (homogeneous) Dirichlet boundary condition at $x_2 = H$ for some artificial height $H>0$~\cite{preston2008integral}. By the method of images, we have that $\widetilde{G}^D$ is given by
\begin{align}
    \label{eq:dirichlet-green-function}
    \widetilde{G}^D(\bx,\by) = &-\frac{1}{4\pi}\log \left((\tau(x_1) - \tau(y_1))^2 + (x_2 - y_2)^2\right) \\
    &+ \frac{1}{4\pi}\log \left((\tau(x_1) - \tau(y_1))^2 + (x_2 - y_2^r)^2 \right) \nonumber
\end{align}
where $y_2^r = 2H - y_2$ is the image point. It is possible to verify that, as in the usual Laplace case, this complexified Dirichlet Green's function decays linearly as $|\by| \to \infty$, which in turn may help define an appropriate functional framework for equations like~\cref{eq:BIE}. 

Secondly, the equivalence between the boundary integral equation~\cref{eq:BIE} and the complex-scaled water-waves problem~\cref{eq:water-waves-system-pml} was not addressed. To briefly indicate the source of difficulty, let us consider a density $\sigma : \Gamma \to \C$ solution of $\cref{eq:BIE}$, and define the function $\tvarphi: \Omega \to \C$ by
\begin{align}
\label{eq:tentative-solution}
\tvarphi(\bx) = \mathcal{D}_{\Gamma}[\sigma](\bx) - \mathcal{S}_{\Gamma}[f - \alpha \tau' \sigma](\bx).
\end{align}
One would like to show that~\cref{eq:tentative-solution} is a solution to the original problem~\cref{eq:water-waves-system-pml}. The fact that~\cref{eq:tentative-solution} satisfies the PDE~\cref{eq:pml-laplace-real} follows directly from the choice of the integral representation. This means that we may write
\begin{align}
\label{eq:tentative-solution-representation-formula}
\tvarphi(\bx) = \mathcal{D}_{\Gamma}[\tvarphi](\bx) - \mathcal{S}_{\Gamma}[\tau' \nabla \tvarphi \cdot \bn](\bx).
\end{align}
Taking the Dirichlet trace of~\cref{eq:tentative-solution} and using the fact that $\sigma$ solves~\cref{eq:BIE} we obtain that $\tvarphi = \sigma$ on $\Gamma$. Subtracting~\cref{eq:tentative-solution-representation-formula} from~\cref{eq:tentative-solution} yields
\begin{align}
\label{eq:zero-potential}
\mathcal{S}_{\Gamma}[\tau' \nabla \tvarphi \cdot \bn - f + \tau'\alpha \tvarphi](\bx) = 0, \quad \bx \in \Omega.
\end{align}
Equality \cref{eq:zero-potential} can be extended to $\Gamma$ to give:
\begin{align}   
\label{eq:tentative-solution-representation-formu}
S_{\Gamma}[\tau' \nabla \tvarphi \cdot \bn - f + \tau'\alpha \tvarphi](\bx) = 0, \quad \bx \in \Gamma.
\end{align}
As a consequence, $\tvarphi$ as given by~\cref{eq:tentative-solution} solves the water-waves system~\cref{eq:water-waves-system-pml} if the single-layer operator $S_\Gamma$ is ``injective" in a certain sense. Due to the current limited understanding of the mapping properties of the single-layer operator in this scenario, we are unable to provide a more precise analysis at this time.

Finally, while the use of a Dirichlet Green's function was motivated above for theoretical reasons, the use of a Green's function satisfying a homogenous Neumann condition at the flat bottom $x_2 = -1$ is more natural from a numerical point of view since it eliminates many degrees of freedom at the bottom (i.e. only the perturbations to the topography need to be discretized when $\Omega$ lies above $x_2=-1$). 
Such a Green's function is easy to construct once again by the method of images:
\begin{align}
    \label{neumann-green-function}
    \widetilde{G}^N(\bx,\by) = &-\frac{1}{4\pi}\log \left((\tau(x_1) - \tau(y_1))^2 + (x_2 - y_2)^2\right) \\
    &- \frac{1}{4\pi}\log \left((\tau(x_1) - \tau(y_1))^2 + (x_2 - y_2^d)^2 \right) \nonumber
\end{align}
where $y_2^d = -2 - y_1$. Unlike the Dirichlet Green's function, however, $\widetilde{G}^N$ still possesses a logarithmic growth at infinity. Because of the aforementioned computational advantage of such a Neumann Green's function, it would be interesting to better understand the properties of the integral equations derived from it. 

\section{Conclusions and future directions}\label{sec:conclusions} 

In this paper, we presented a novel boundary integral formulation for solving the two-dimensional finite depth water waves problem. We showed that the underlying truncation methodology leads to exponentially small truncation errors, and it is applicable to many challenging problems including scattering in the presence of a step topography. We believe the ability to compute resonances (both real and complex) in the form of a linear eigenvalue problem, while enjoying the dimensional reduction of boundary integral equations, is a very interesting feature of the method worthy of further exploration. Future directions of research include extensions to infinite depth and the three-dimensional problem. 

\appendix

\section{Complex-scaled fundamental solution}\label{sec:fundamental-solution}
Let $x_1\in\R\mapsto \tau(x_1)\in\C$ be as in section \ref{sec:complex-stretching} a   continuous and piecewise continuously differentiable function satisfying \eqref{eq:tau-condition-real}. 
 In this section, we provide a direct proof that $$\widetilde{G}(\bx,\by)=G(\btau(\bx),\btau(\by))\mbox{ where }\btau(\bx)=(\tau(x_1),x_2)$$  is the free-space fundamental solution of the complex-scaled Laplace equation,  i.e. that $\widetilde{G}$
 satisfies 
 \begin{align}
 \label{eq:identityGreen}
  \int_{\R^2} \nabla_\by\widetilde{G}(\bx,\by) \cdot   A(\by) \nabla v(\by)   \de \by = v(\bx),\quad \forall \bx\in \R^2,
 \end{align}
 for all test functions $v \in \mathcal{D}(\R^2)$, where $A$ is given by \eqref{eq:water-waves-system-pml}.
 
By continuity arguments, it suffices to prove \eqref{eq:identityGreen} for points $\bx$ where $\btau$ is differentiable. Note that  $\widetilde{G}(\bx,\by)$ is well-defined for all $\bx$, $\by$ in $\R^2$, as soon as $\bx\neq\by$. Indeed, suppose that 
$$(\tau(x_1)-\tau(y_1))^2+(x_2-y_2)^2=0,$$
then 
$$\real (\tau(x_1)-\tau(y_1))=0,$$
which implies by \eqref{eq:tau-condition-real} that $x_1=y_1$, and finally that $\bx=\by$. 
Now by the
 properties of the Lebesgue integral and the symmetry of $A(\by)$, it holds that, by letting $I(\bx)$ denote the left hand-side of~\cref{eq:identityGreen},
 $$I(\bx)=\lim_{\varepsilon\to 0}\int_{|\by-\bx|>\varepsilon}  A(\by)\nabla_\by\widetilde{G}(\bx,\by) \cdot  \nabla v(\by)    \de \by.$$
 Integrating by parts and using that 
 $$\nabla_{\by} \cdot( A({\by}) \nabla_\by\widetilde{G}(\bx,\by))=0\mbox{ if }\by\neq \bx,$$
 we obtain that
 $$I(\bx)=\lim_{\varepsilon\to 0}\frac{1}{\varepsilon}\int_{|\by-\bx|=\varepsilon}A({\by}) \nabla_\by\widetilde{G}(\bx,\by)\cdot (\bx-\by)v(\by)    \de s_\by.$$
 On the other hand, the chain rule and the definition of the function $G$ lead to 
 $$\nabla_\by\widetilde{G}(\bx,\by)=\frac{1}{2\pi }\frac{1}{(\tau(x_1)-\tau(y_1))^2+(x_2-y_2)^2}\begin{bmatrix}
  \tau'(y_1)(\tau(x_1)-\tau(y_1))  \\
    x_2-y_2 
 \end{bmatrix},$$
 so that by definition of $A(\by)$:
  $$I(\bx)=\lim_{\varepsilon\to 0}\frac{1}{2\pi\varepsilon}\int_{|\by-\bx|=\varepsilon}\!\!\!\tau'(y_1)\frac{ \tau'(y_1)^{-1}(\tau(x_1)-\tau(y_1))(x_1-y_1) +(x_2-y_2)^2}{(\tau(x_1)-\tau(y_1))^2+(x_2-y_2)^2}v(\by)  \de s_\by.$$
  Now parametrizing the circle by $\by = \bx + \epsilon\br(\theta)$, with
 $\br(\theta) = (\cos(\theta),\sin(\theta))$ and using straightforward Taylor expansions, this leads to 
 $$I(\bx)=\frac{v(\bx)}{2\pi}\int_0^{2\pi} \frac{\tau'(x_1)}{\tau'(x_1)^2\cos(\theta)^2+\sin(\theta)^2} \de \theta.$$
We conclude by using the following classical result of complex analysis: 
$$\int_0^{2\pi} \frac{1}{\alpha^2\cos(\theta)^2+\sin(\theta)^2} \de \theta=\frac{2\pi}{\alpha}$$
which holds for all $\alpha\in\C$ such that $\real(\alpha)>0$. The result can be proved first for real $\alpha$, rewriting the identity as
$$\imag \left(\int_{\mathcal C}\frac{dz}{z}\right)=2\pi$$
where ${\mathcal C}$ is an ellpise defined by  ${\mathcal C}=\{\alpha\cos(\theta)+i\sin(\theta); 0\leq\theta\leq 2\pi\}$. Then one concludes that the identity still holds for complex $\alpha$ by analyticity.

\bibliographystyle{siamplain}
\bibliography{references}

\end{document}